\documentclass[12pt,thmsa]{article}%
\usepackage{amsfonts}
\usepackage{sw20bams}
\usepackage{amsmath}
\usepackage{amssymb}
\usepackage{graphicx}%
\providecommand{\U}[1]{\protect\rule{.1in}{.1in}}
\begin{document}

\author{Steven Finch}
\title{Quartic and Octic Characters Modulo $n$}
\date{March 25, 2016}
\maketitle

\begin{abstract}
The average number of primitive quadratic Dirichlet characters of modulus $n$
tends to a constant as $n\rightarrow\infty$. The same is true for primitive
cubic characters. It is therefore surprising that, as $n\rightarrow\infty$,
the average number of primitive quartic characters of modulus $n$ grows with
$\ln(n)$, and that the average number of primitive octic characters of modulus
$n$ grows with $\ln(n)^{2}$. Leading coefficients in the asymptotic
expressions are also computed.

\end{abstract}

\footnotetext{Copyright \copyright \ 2009 by Steven R. Finch. All rights
reserved.}Let $\mathbb{Z}_{n}^{\ast}$ denote the group (under multiplication
modulo $n$) of integers relatively prime to $n$, and let $\mathbb{C}^{\ast}$
denote the group (under ordinary multiplication) of nonzero complex numbers.
We wish to count homomorphisms $\chi:\mathbb{Z}_{n}^{\ast}\rightarrow
\mathbb{C}^{\ast}$ satisfying certain requirements. A Dirichlet character
$\chi$ is quadratic if $\chi(k)^{2}=1$ for every $k$ in $\mathbb{Z}_{n}^{\ast
}$. It is well-known that
\[%
\begin{array}
[c]{l}%
\text{\#\ quadratic Dirichlet characters}\\
\text{of modulus }\leq N
\end{array}
=%
{\displaystyle\sum\limits_{n\leq N}}
a(n)\sim\frac{6}{\pi^{2}}N\ln(N)
\]
as $N\rightarrow\infty$, where $a(n)$ is multiplicative with
\[%
\begin{array}
[c]{ccc}%
a(2^{r})=\left\{
\begin{array}
[c]{lll}%
1 &  & \text{if }r=1,\\
2 &  & \text{if }r=2,\\
4 &  & \text{if }r\geq3,
\end{array}
\right.  &  & a(p^{r})=2
\end{array}
\]
for prime $p\geq3$ and $r\geq1$; and
\[%
\begin{array}
[c]{l}%
\text{\#\ primitive quadratic Dirichlet }\\
\text{characters of modulus }\leq N
\end{array}
=%
{\displaystyle\sum\limits_{n\leq N}}
b(n)\sim\frac{6}{\pi^{2}}N
\]
where $b(n)$ is multiplicative with
\[%
\begin{array}
[c]{ccc}%
b(2^{r})=\left\{
\begin{array}
[c]{lll}%
0 &  & \text{if }r=1,\\
1 &  & \text{if }r=2,\\
2 &  & \text{if }r=3,\\
0 &  & \text{if }r\geq4,
\end{array}
\right.  &  & b(p^{r})=\left\{
\begin{array}
[c]{lll}%
1 &  & \text{if }r=1,\\
0 &  & \text{if }r\geq2
\end{array}
\right.
\end{array}
\]
for prime $p\geq3$. Is it surprising that the same constant $6/\pi^{2}$
appears in connection with averages of both $a(n)$ and $b(n)$? No. The
underlying structure linking these two sequences is the M\"{o}bius inversion
formula \cite{Wlf}:
\[%
\begin{array}
[c]{ccc}%
{\displaystyle\sum\limits_{n=1}^{\infty}}
\dfrac{a(n)}{n^{s}}=\zeta(s)%
{\displaystyle\sum\limits_{n=1}^{\infty}}
\dfrac{b(n)}{n^{s}} &  & \text{for }\operatorname*{Re}(s)>c\text{ for some
}c<1.
\end{array}
\]
This can be proved directly here. On the one hand \cite{Sec0},
\begin{align*}%
{\displaystyle\sum\limits_{n=1}^{\infty}}
\frac{a(n)}{n^{s}}  &  =\left(  1+\frac{1}{2^{s}}+\frac{2}{2^{2s}}+%
{\displaystyle\sum\limits_{r=3}^{\infty}}
\frac{4}{2^{rs}}\right)
{\displaystyle\prod\limits_{p>2}}
\left(  1+%
{\displaystyle\sum\limits_{r=1}^{\infty}}
\frac{2}{p^{rs}}\right) \\
\  &  =\left(  1+\frac{1}{2^{s}}+\frac{2}{2^{2s}}+\frac{4}{2^{s}-1}-\frac
{4}{2^{s}}-\frac{4}{2^{2s}}\right)
{\displaystyle\prod\limits_{p>2}}
\left(  1+\frac{2}{p^{s}-1}\right) \\
\  &  =\left(  1-\frac{3}{2^{s}}-\frac{2}{2^{2s}}+\frac{4}{2^{s}-1}\right)
{\displaystyle\prod\limits_{p>2}}
\left(  1-\frac{1}{p^{s}}\right)  ^{-2}\left(  1-\frac{1}{p^{2s}}\right) \\
\  &  =\left(  1-\frac{3}{2^{s}}-\frac{2}{2^{2s}}+\frac{4}{2^{s}-1}\right)
\left(  1-\frac{1}{2^{s}}\right)  ^{2}\zeta(s)^{2}\left(  1-\frac{1}{2^{2s}%
}\right)  ^{-1}\zeta(2s)^{-1}\\
\  &  =\left(  1+\frac{1}{2^{2s}}+\frac{2}{2^{3s}}\right)  \left(  1-\frac
{1}{2^{s}}\right)  \zeta(s)^{2}\left(  1-\frac{1}{2^{2s}}\right)  ^{-1}%
\zeta(2s)^{-1};
\end{align*}
on the other hand,
\begin{align*}%
{\displaystyle\sum\limits_{n=1}^{\infty}}
\frac{b(n)}{n^{s}}  &  =\left(  1+\frac{1}{2^{2s}}+\frac{2}{2^{3s}}\right)
{\displaystyle\prod\limits_{p>2}}
\left(  1+\frac{1}{p^{s}}\right) \\
\  &  =\left(  1+\frac{1}{2^{2s}}+\frac{2}{2^{3s}}\right)
{\displaystyle\prod\limits_{p>2}}
\left(  1-\frac{1}{p^{s}}\right)  ^{-1}\left(  1-\frac{1}{p^{2s}}\right) \\
\  &  =\left(  1+\frac{1}{2^{2s}}+\frac{2}{2^{3s}}\right)  \left(  1-\frac
{1}{2^{s}}\right)  \zeta(s)\left(  1-\frac{1}{2^{2s}}\right)  ^{-1}%
\zeta(2s)^{-1}.
\end{align*}
The leading coefficient $6/\pi^{2}$ arises because
\[
\left.  \left(  1+\frac{1}{2^{2s}}+\frac{2}{2^{3s}}\right)  \left(  1-\frac
{1}{2^{s}}\right)  \left(  1-\frac{1}{2^{2s}}\right)  ^{-1}\zeta
(2s)^{-1}\right\vert _{s=1}=\frac{1}{\zeta(2)}=\frac{6}{\pi^{2}}.
\]
The M\"{o}bius inversion formula is valid, in fact, for arbitrary
$\ell^{\text{th}}$-order Dirichlet characters. From the preceding case
$\ell=2$, we observe the evaluation of $\sum_{n=1}^{\infty}b(n)n^{-s}$ to be
slightly less complicated than that for $\sum_{n=1}^{\infty}a(n)n^{-s}$.
Hence, to simplify calculations, our focus will be on $b(n)$ for the cases
$\ell=3,4,8$.

\section{Cubic Characters}

\subsection{General}

A Dirichlet character $\chi$ is cubic if $\chi(k)^{3}=1$ for every $k$ in
$\mathbb{Z}_{n}^{*}$. We have \cite{Sec1}
\[%
\begin{array}
[c]{ccccc}%
a(2^{r})=1, &  & a(3^{r})=\left\{
\begin{array}
[c]{lll}%
1 &  & \text{if }r=1,\\
3 &  & \text{if }r\geq2,
\end{array}
\right.  &  & a(p^{r})=\left\{
\begin{array}
[c]{lll}%
1 &  & \text{if }p\equiv2\operatorname*{mod}3,\\
3 &  & \text{if }p\equiv1\operatorname*{mod}3
\end{array}
\right.
\end{array}
\]
for prime $p\geq5$ and $r\geq1$. The asymptotics for $%
{\textstyle\sum\nolimits_{n\leq N}}
a(n)$, as well as the coefficient, were studied in \cite{FS} (via a different
proof than the following).

\subsection{Primitive}

We have%

\[%
\begin{array}
[c]{ccccc}%
b(2^{r})=0, &  & b(3^{r})=\left\{
\begin{array}
[c]{lll}%
2 &  & \text{if }r=2,\\
0 &  & \text{otherwise,}%
\end{array}
\right.  &  & b(p^{r})=\left\{
\begin{array}
[c]{lll}%
2 &  & \text{if }r=1\text{ \& }p\equiv1\operatorname*{mod}3,\\
0 &  & \text{otherwise}%
\end{array}
\right.
\end{array}
\]
for prime $p\geq5$ and $r\geq1$. The asymptotics for $%
{\textstyle\sum\nolimits_{n\leq N}}
b(n)$, as well as the coefficient, were studied in \cite{DFK}. We obtain
\begin{align*}%
{\displaystyle\sum\limits_{n=1}^{\infty}}
\frac{b(n)}{n^{s}}  &  =\left(  1+\frac2{3^{2s}}\right)
{\displaystyle\prod\limits_{p\equiv1\operatorname*{mod}3}}
\left(  1+\frac2{p^{s}}\right) \\
\  &  =\left(  1+\frac2{3^{2s}}\right)
{\displaystyle\prod\limits_{p\equiv1\operatorname*{mod}3}}
\left(  1-\frac1{p^{s}}\right)  ^{-2}\left(  1-\frac1{p^{2s}}\right)  \left(
1-\frac2{p^{s}(p^{s}+1)}\right)
\end{align*}
and it is known from section [\ref{Euler3}] that
\begin{equation}%
{\displaystyle\prod\limits_{p\equiv1\operatorname*{mod}3}}
\left(  1-\frac1{p^{s}}\right)  ^{-2}\sim\zeta(s)\frac{\sqrt{3}}{2\pi}%
{\displaystyle\prod\limits_{p\equiv1\operatorname*{mod}3}}
\left(  1-\frac1{p^{2}}\right)  ^{-1} \label{E1}%
\end{equation}
as $s\rightarrow1$; after cancellation, the coefficient becomes
\[
\frac{11}9\frac{\sqrt{3}}{2\pi}%
{\displaystyle\prod\limits_{p\equiv1\operatorname*{mod}3}}
\left(  1-\frac2{p(p+1)}\right)  =0.3170565167922841205670156....
\]
Just as in the case of primitive quadratic characters, the function $\zeta(s)$
appears with exponent $1$. Hence by the Selberg-Delange method \cite{FS}, the
average number of primitive cubic characters of modulus $n$ tends to a
constant $0.317...$ as $n\rightarrow\infty$.

\section{Quartic Characters}

\subsection{General}

A Dirichlet character $\chi$ is quartic (biquadratic) if $\chi(k)^{4}=1$ for
every $k$ in $\mathbb{Z}_{n}^{*}$. We have \cite{Sec2}
\[%
\begin{array}
[c]{ccc}%
a(2^{r})=\left\{
\begin{array}
[c]{lll}%
1 &  & \text{if }r=1,\\
2 &  & \text{if }r=2,\\
4 &  & \text{if }r=3,\\
8 &  & \text{if }r\geq4,
\end{array}
\right.  &  & a(p^{r})=\left\{
\begin{array}
[c]{lll}%
2 &  & \text{if }p\equiv3\operatorname*{mod}4,\\
4 &  & \text{if }p\equiv1\operatorname*{mod}4
\end{array}
\right.
\end{array}
\]
for prime $p\geq3$ and $r\geq1$.

\subsection{Primitive}

We have%

\[%
\begin{array}
[c]{ccc}%
b(2^{r})=\left\{
\begin{array}
[c]{lll}%
1 &  & \text{if }r=2,\\
2 &  & \text{if }r=3,\\
4 &  & \text{if }r=4,\\
0 &  & \text{otherwise,}%
\end{array}
\right.  &  & b(p^{r})=\left\{
\begin{array}
[c]{lll}%
1 &  & \text{if }r=1\text{ \& }p\equiv3\operatorname*{mod}4,\\
3 &  & \text{if }r=1\text{ \& }p\equiv1\operatorname*{mod}4\\
0 &  & \text{otherwise}%
\end{array}
\right.
\end{array}
\]
for prime $p\geq3$ and $r\geq1$. We obtain
\begin{align*}%
{\displaystyle\sum\limits_{n=1}^{\infty}}
\frac{b(n)}{n^{s}}  &  =\left(  1+\frac1{2^{2s}}+\frac2{2^{3s}}+\frac4{2^{4s}%
}\right)
{\displaystyle\prod\limits_{p\equiv3\operatorname*{mod}4}}
\left(  1+\frac1{p^{s}}\right)  \cdot%
{\displaystyle\prod\limits_{p\equiv1\operatorname*{mod}4}}
\left(  1+\frac3{p^{s}}\right) \\
\  &  =\left(  1+\frac1{2^{2s}}+\frac2{2^{3s}}+\frac4{2^{4s}}\right)
{\displaystyle\prod\limits_{p\equiv3\operatorname*{mod}4}}
\left(  1-\frac1{p^{s}}\right)  ^{-1}\left(  1-\frac1{p^{2s}}\right)  \cdot\\
&  \ \ \cdot%
{\displaystyle\prod\limits_{p\equiv1\operatorname*{mod}4}}
\left(  1-\frac1{p^{s}}\right)  ^{-3}\left(  1-\frac1{p^{2s}}\right)  \left(
1-\frac{5p^{s}-3}{p^{2s}(p^{s}+1)}\right) \\
\  &  =\left(  1+\frac1{2^{2s}}+\frac2{2^{3s}}+\frac4{2^{4s}}\right)  \left(
1-\frac1{2^{s}}\right)  \zeta(s)\left(  1-\frac1{2^{2s}}\right)  ^{-1}%
\zeta(2s)^{-1}\cdot\\
&  \ \ \cdot%
{\displaystyle\prod\limits_{p\equiv1\operatorname*{mod}4}}
\left(  1-\frac1{p^{s}}\right)  ^{-2}\left(  1-\frac{5p^{s}-3}{p^{2s}%
(p^{s}+1)}\right)
\end{align*}
and, from section [\ref{Euler4}],
\begin{equation}%
{\displaystyle\prod\limits_{p\equiv1\operatorname*{mod}4}}
\left(  1-\frac1{p^{s}}\right)  ^{-2}\sim\zeta(s)\frac1\pi%
{\displaystyle\prod\limits_{p\equiv1\operatorname*{mod}4}}
\left(  1-\frac1{p^{2}}\right)  ^{-1}\sim\zeta(s)\frac\pi{16K^{2}} \label{E2}%
\end{equation}
as $s\rightarrow1$; thus the coefficient becomes
\[
\frac7\pi\frac1{16K^{2}}%
{\displaystyle\prod\limits_{p\equiv1\operatorname*{mod}4}}
\left(  1-\frac{5p-3}{p^{2}(p+1)}\right)  =0.1908767211685284480112237....
\]
where $K$ is the Landau-Ramanujan constant \cite{Fnch}. Unlike the preceding
cases, the function $\zeta(s)$ appears with exponent $2$. Hence by the
Selberg-Delange method \cite{FS}, the average number of primitive quartic
characters of modulus $n$ is asymptotically $(0.190...)\ln(n)$ as
$n\rightarrow\infty$. Is it surprising that the quartic case differs so
dramatically from both the quadratic and cubic cases? We believe yes. There is
no \textit{a priori} reason for quartic characters to outnumber
quadratic/cubic characters in such a manner.

\section{Octic Characters}

\subsection{General}

A Dirichlet character $\chi$ is octic if $\chi(k)^{8}=1$ for every $k$ in
$\mathbb{Z}_{n}^{*}$. We have
\[%
\begin{array}
[c]{ccc}%
a(2^{r})=\left\{
\begin{array}
[c]{lll}%
1 &  & \text{if }r=1,\\
2 &  & \text{if }r=2,\\
4 &  & \text{if }r=3,\\
8 &  & \text{if }r=4,\\
16 &  & \text{if }r\geq5,
\end{array}
\right.  &  & a(p^{r})=\left\{
\begin{array}
[c]{lll}%
2 &  & \text{if }p\equiv3\operatorname*{mod}4,\\
4 &  & \text{if }p\equiv5\operatorname*{mod}8,\\
8 &  & \text{if }p\equiv1\operatorname*{mod}8
\end{array}
\right.
\end{array}
\]
for prime $p\geq3$ and $r\geq1$.

\subsection{Primitive}

We have%

\[%
\begin{array}
[c]{ccc}%
b(2^{r})=\left\{
\begin{array}
[c]{lll}%
1 &  & \text{if }r=2,\\
2 &  & \text{if }r=3,\\
4 &  & \text{if }r=4,\\
8 &  & \text{if }r=5,\\
0 &  & \text{otherwise,}%
\end{array}
\right.  &  & b(p^{r})=\left\{
\begin{array}
[c]{lll}%
1 &  & \text{if }r=1\text{ \& }p\equiv3\operatorname*{mod}4,\\
3 &  & \text{if }r=1\text{ \& }p\equiv5\operatorname*{mod}8,\\
7 &  & \text{if }r=1\text{ \& }p\equiv1\operatorname*{mod}8,\\
0 &  & \text{otherwise}%
\end{array}
\right.
\end{array}
\]
for prime $p\geq3$ and $r\geq1$. We obtain
\begin{align*}%
{\textstyle\sum\limits_{n=1}^{\infty}}
\tfrac{b(n)}{n^{s}}  &  =\left(  1+\tfrac1{2^{2s}}+\tfrac2{2^{3s}}%
+\tfrac4{2^{4s}}+\tfrac8{2^{5s}}\right)
{\textstyle\prod\limits_{\substack{p\equiv3 \\\operatorname*{mod}4 }}}
\left(  1+\tfrac1{p^{s}}\right)  \cdot%
{\textstyle\prod\limits_{\substack{p\equiv5 \\\operatorname*{mod}8 }}}
\left(  1+\tfrac3{p^{s}}\right)  \cdot%
{\textstyle\prod\limits_{\substack{p\equiv1 \\\operatorname*{mod}8 }}}
\left(  1+\tfrac7{p^{s}}\right) \\
\  &  =\left(  1+\tfrac1{2^{2s}}+\tfrac2{2^{3s}}+\tfrac4{2^{4s}}%
+\tfrac8{2^{5s}}\right)
{\textstyle\prod\limits_{\substack{p\equiv3 \\\operatorname*{mod}4 }}}
\left(  1-\tfrac1{p^{s}}\right)  ^{-1}\left(  1-\tfrac1{p^{2s}}\right)
\cdot\\
&  \ \cdot%
{\textstyle\prod\limits_{\substack{p\equiv5 \\\operatorname*{mod}8 }}}
\left(  1-\tfrac1{p^{s}}\right)  ^{-3}\left(  1-\tfrac1{p^{2s}}\right)
\left(  1-\tfrac{5p^{s}-3}{p^{2s}(p^{s}+1)}\right) \\
&  \cdot%
{\textstyle\prod\limits_{\substack{p\equiv1 \\\operatorname*{mod}8 }}}
\left(  1-\tfrac1{p^{s}}\right)  ^{-7}\left(  1-\tfrac1{p^{2s}}\right)
\left(  1-\tfrac{27p^{5s}-85p^{4s}+125p^{3s}-99p^{2s}+41p^{s}-7}{p^{6s}%
(p^{s}+1)}\right) \\
\  &  =\left(  1+\tfrac1{2^{2s}}+\tfrac2{2^{3s}}+\tfrac4{2^{4s}}%
+\tfrac8{2^{5s}}\right)  \left(  1-\tfrac1{2^{s}}\right)  \zeta(s)\left(
1-\tfrac1{2^{2s}}\right)  ^{-1}\zeta(2s)^{-1}\cdot\\
&  \ \cdot%
{\textstyle\prod\limits_{\substack{p\equiv5 \\\operatorname*{mod}8 }}}
\left(  1-\tfrac1{p^{s}}\right)  ^{-2}\left(  1-\tfrac{5p^{s}-3}{p^{2s}%
(p^{s}+1)}\right)  \cdot\\
&  \cdot%
{\textstyle\prod\limits_{\substack{p\equiv1 \\\operatorname*{mod}8 }}}
\left(  1-\tfrac1{p^{s}}\right)  ^{-6}\left(  1-\tfrac{27p^{5s}-85p^{4s}%
+125p^{3s}-99p^{2s}+41p^{s}-7}{p^{6s}(p^{s}+1)}\right)
\end{align*}
and, from section [\ref{Euler8}],
\begin{equation}%
{\displaystyle\prod\limits_{p\equiv5\operatorname*{mod}8}}
\left(  1-\frac1{p^{s}}\right)  ^{-4}\sim\zeta(s)\frac1{2\ln\left(  1+\sqrt
{2}\right)  }%
{\displaystyle\prod\limits_{p\equiv5\operatorname*{mod}8}}
\left(  1-\frac1{p^{2}}\right)  ^{-2}, \label{E3}%
\end{equation}
\begin{equation}%
{\displaystyle\prod\limits_{p\equiv1\operatorname*{mod}8}}
\left(  1-\frac1{p^{s}}\right)  ^{-4}\sim\zeta(s)\frac{2\ln\left(  1+\sqrt
{2}\right)  }{\pi^{2}}%
{\displaystyle\prod\limits_{p\equiv1\operatorname*{mod}8}}
\left(  1-\frac1{p^{2}}\right)  ^{-2} \label{E4}%
\end{equation}
as $s\rightarrow1$. An expression for the coefficient becomes clear. More
importantly, the function $\zeta(s)$ appears with exponent $1+1/2+3/2=3$.
Hence by the Selberg-Delange method \cite{FS}, the average number of primitive
octic characters of modulus $n$ has growth rate $\approx\ln(n)^{2}$ as
$n\rightarrow\infty$.

\section{Euler Product Residues}

Formulas (\ref{E1}), (\ref{E2}), (\ref{E4}) await proof, while the truth of
(\ref{E3}) depends on both (\ref{E2}) and (\ref{E4}). Our approach uses the
seemingly-unrelated method of Shanks \&\ Schmid \cite{SS} for computing
various generalized Landau-Ramanujan constants $\kappa_{m}$.

Fix an integer $m\neq0$. Define
\[
L_{d}(s)=%
{\displaystyle\sum\limits_{n=1}^{\infty}}
\genfrac{(}{)}{}{}{d}{n}%
\frac1{n^{s}}
\]
where $d=-m$ if $4\,|\,m$ and $d=-4m$ otherwise; $(\cdot/\cdot)$ is the
Kronecker-Jacobi-Legendre symbol. Define also
\[
\Lambda_{m}(s)=%
{\displaystyle\sum\limits_{n=1}^{\infty}}
\frac{f(n)}{n^{s}}
\]
where $f(n)=1$ if there exist integers $x,y$ such that $n=x^{2}+m\,y^{2}$ and
$f(n)=0$ otherwise. Let
\[
\kappa_{m}=\delta_{m}\sqrt{\frac{L_{d}(1)}\pi\frac{2|m|}{\varphi\left(
2|m|\right)  }}%
{\displaystyle\prod\limits_{\substack{\text{odd }p: \\\binom{-m}{p}=-1 }}}
\left(  1-\frac1{p^{2}}\right)  ^{-1/2}
\]
where the rational number $\delta_{m}$ is unspecified for the moment, and
$\varphi$ denotes the Euler totient function. It turns out that $\sum_{n\leq
N}f(n)\sim\kappa_{m}N/\sqrt{\ln(N)}$ as $N\rightarrow\infty$, although this
fact is not directly relevant to our purposes.

\subsection{\label{Euler3}Case when $\ell=3$}

Let $m=3$ and $\delta_{m}=2/3$. It follows that
\begin{align*}
\kappa_{3}  &  =\frac23\sqrt{L_{-12}(1)\frac1\pi\frac6{\varphi\left(
6\right)  }}%
{\displaystyle\prod\limits_{\substack{\text{odd }p: \\p\equiv
2\operatorname*{mod}3 }}}
\left(  1-\frac1{p^{2}}\right)  ^{-1/2}\\
\  &  =\frac23\sqrt{\frac\pi{2\sqrt{3}}\frac1\pi\frac62}%
{\displaystyle\prod\limits_{p\equiv2\operatorname*{mod}3}}
\left(  1-\frac1{p^{2}}\right)  ^{-1/2}\cdot\left(  1-\frac1{2^{2}}\right)
^{1/2}\\
\  &  =\frac1{\sqrt{2\sqrt{3}}}%
{\displaystyle\prod\limits_{p\equiv2\operatorname*{mod}3}}
\left(  1-\frac1{p^{2}}\right)  ^{-1/2}=2^{1/2}3^{-7/4}\pi%
{\displaystyle\prod\limits_{p\equiv1\operatorname*{mod}3}}
\left(  1-\frac1{p^{2}}\right)  ^{1/2}.
\end{align*}
Now,
\[
\Lambda_{3}(s)=\left(  1-3^{-s}\right)  ^{-1}%
{\displaystyle\prod\limits_{\substack{p\equiv1 \\\operatorname*{mod}3 }}}
\left(  1-p^{-s}\right)  ^{-1}\cdot%
{\displaystyle\prod\limits_{\substack{p\equiv2 \\\operatorname*{mod}3 }}}
\left(  1-p^{-2s}\right)  ^{-1}
\]
by definition and
\[
\Lambda_{3}(s)^{2}=\zeta(s)L_{-3}(s)\left(  1-3^{-s}\right)  ^{-1}%
{\displaystyle\prod\limits_{\substack{p\equiv2 \\\operatorname*{mod}3 }}}
\left(  1-p^{-2s}\right)  ^{-1}
\]
by elementary considerations \cite{MR}. Dividing the second by the first, we
have
\[
\Lambda_{3}(s)=\zeta(s)L_{-3}(s)%
{\displaystyle\prod\limits_{\substack{p\equiv1 \\\operatorname*{mod}3 }}}
\left(  1-p^{-s}\right)
\]
and thus
\[
\underset{
\begin{array}
[c]{c}%
\downarrow\\
\sqrt{\pi}\kappa_{3}%
\end{array}
}{\underbrace{\Lambda_{3}(s)(s-1)^{1/2}}}\,=\,\underset{
\begin{array}
[c]{c}%
\downarrow\\
1
\end{array}
}{\underbrace{\zeta(s)(s-1)}\,}\underset{
\begin{array}
[c]{c}%
\downarrow\\
\tfrac\pi{3\sqrt{3}}%
\end{array}
}{\underbrace{L_{-3}(s)}}%
{\displaystyle\prod\limits_{\substack{p\equiv1 \\\operatorname*{mod}3 }}}
\left(  1-p^{-s}\right)  \cdot(s-1)^{-1/2}
\]
as $s\rightarrow1$. Hence formula (\ref{E1}) is true.

\subsection{\label{Euler4}Case when $\ell=4$}

Let $m=1$ and $\delta_{m}=1$. It follows that
\begin{align*}
\kappa_{1}  &  =1\sqrt{L_{-4}(1)\frac1\pi\frac2{\varphi\left(  2\right)  }}%
{\displaystyle\prod\limits_{\substack{\text{odd }p: \\p\equiv
3\operatorname*{mod}4 }}}
\left(  1-\frac1{p^{2}}\right)  ^{-1/2}\\
\  &  =\sqrt{\frac\pi4\frac1\pi\frac21}%
{\displaystyle\prod\limits_{p\equiv3\operatorname*{mod}4}}
\left(  1-\frac1{p^{2}}\right)  ^{-1/2}\\
\  &  =\frac1{\sqrt{2}}%
{\displaystyle\prod\limits_{p\equiv3\operatorname*{mod}4}}
\left(  1-\frac1{p^{2}}\right)  ^{-1/2}=K=\frac\pi4%
{\displaystyle\prod\limits_{p\equiv1\operatorname*{mod}4}}
\left(  1-\frac1{p^{2}}\right)  ^{1/2}.
\end{align*}
Now,
\[
\Lambda_{1}(s)=\left(  1-2^{-s}\right)  ^{-1}%
{\displaystyle\prod\limits_{\substack{p\equiv1 \\\operatorname*{mod}4 }}}
\left(  1-p^{-s}\right)  ^{-1}\cdot%
{\displaystyle\prod\limits_{\substack{p\equiv3 \\\operatorname*{mod}4 }}}
\left(  1-p^{-2s}\right)  ^{-1}
\]
by definition and
\[
\Lambda_{1}(s)^{2}=\zeta(s)L_{-4}(s)\left(  1-2^{-s}\right)  ^{-1}%
{\displaystyle\prod\limits_{\substack{p\equiv3 \\\operatorname*{mod}4 }}}
\left(  1-p^{-2s}\right)  ^{-1}
\]
by elementary considerations \cite{Hrdy}. Dividing the second by the first, we
obtain
\[
\underset{
\begin{array}
[c]{c}%
\downarrow\\
\sqrt{\pi}\kappa_{1}%
\end{array}
}{\underbrace{\Lambda_{1}(s)(s-1)^{1/2}}}\,=\,\underset{
\begin{array}
[c]{c}%
\downarrow\\
1
\end{array}
}{\underbrace{\zeta(s)(s-1)}}\,\underset{
\begin{array}
[c]{c}%
\downarrow\\
\tfrac\pi4
\end{array}
}{\underbrace{L_{-4}(s)}}%
{\displaystyle\prod\limits_{\substack{p\equiv1 \\\operatorname*{mod}4 }}}
\left(  1-p^{-s}\right)  \cdot(s-1)^{-1/2}
\]
as $s\rightarrow1$. Hence formula (\ref{E2}) is true.

\subsection{\label{Euler8}Case when $\ell=8$}

The argument here seems fairly roundabout and we wonder if a simpler approach
is possible. Two values of $m$ require consideration here. First, let $m=2$
and $\delta_{m}=1$. It follows that
\begin{align*}
\kappa_{2}  &  =1\sqrt{L_{-8}(1)\frac1\pi\frac4{\varphi\left(  4\right)  }}%
{\displaystyle\prod\limits_{\substack{\text{odd }p: \\p\equiv
5,7\operatorname*{mod}8 }}}
\left(  1-\frac1{p^{2}}\right)  ^{-1/2}\\
\  &  =\sqrt{\frac\pi{2\sqrt{2}}\frac1\pi\frac42}%
{\displaystyle\prod\limits_{p\equiv5,7\operatorname*{mod}8}}
\left(  1-\frac1{p^{2}}\right)  ^{-1/2}\\
\  &  =\frac1{\sqrt[4]{2}}%
{\displaystyle\prod\limits_{p\equiv5,7\operatorname*{mod}8}}
\left(  1-\frac1{p^{2}}\right)  ^{-1/2}=2^{-7/4}\pi%
{\displaystyle\prod\limits_{p\equiv1,3\operatorname*{mod}8}}
\left(  1-\frac1{p^{2}}\right)  ^{1/2}.
\end{align*}
As before,
\[
\Lambda_{2}(s)=\left(  1-2^{-s}\right)  ^{-1}%
{\displaystyle\prod\limits_{\substack{p\equiv1,3 \\\operatorname*{mod}8 }}}
\left(  1-p^{-s}\right)  ^{-1}\cdot%
{\displaystyle\prod\limits_{\substack{p\equiv5,7 \\\operatorname*{mod}8 }}}
\left(  1-p^{-2s}\right)  ^{-1},
\]
\[
\Lambda_{2}(s)^{2}=\zeta(s)L_{-8}(s)\left(  1-2^{-s}\right)  ^{-1}%
{\displaystyle\prod\limits_{\substack{p\equiv5,7 \\\operatorname*{mod}8 }}}
\left(  1-p^{-2s}\right)  ^{-1}
\]
and, upon division, we obtain
\[
\underset{
\begin{array}
[c]{c}%
\downarrow\\
\sqrt{\pi}\kappa_{2}%
\end{array}
}{\underbrace{\Lambda_{2}(s)(s-1)^{1/2}}}\,=\,\underset{
\begin{array}
[c]{c}%
\downarrow\\
1
\end{array}
}{\underbrace{\zeta(s)(s-1)}}\,\underset{
\begin{array}
[c]{c}%
\downarrow\\
\tfrac\pi{2\sqrt{2}}%
\end{array}
}{\underbrace{L_{-8}(s)}}%
{\displaystyle\prod\limits_{\substack{p\equiv1,3 \\\operatorname*{mod}8 }}}
\left(  1-p^{-s}\right)  \cdot(s-1)^{-1/2}.
\]

\newpage\ Second, let $m=-2$ and $\delta_{m}=1$. It follows that
\begin{align*}
\kappa_{-2}  &  =1\sqrt{L_{8}(1)\frac1\pi\frac4{\varphi\left(  4\right)  }}%
{\displaystyle\prod\limits_{\substack{\text{odd }p: \\p\equiv
3,5\operatorname*{mod}8 }}}
\left(  1-\frac1{p^{2}}\right)  ^{-1/2}\\
\  &  =\sqrt{\frac{\sqrt{2}\ln\left(  1+\sqrt{2}\right)  }\pi}%
{\displaystyle\prod\limits_{p\equiv3,5\operatorname*{mod}8}}
\left(  1-\frac1{p^{2}}\right)  ^{-1/2}\\
\  &  =2^{-5/4}\ln\left(  1+\sqrt{2}\right)  ^{1/2}\pi^{1/2}%
{\displaystyle\prod\limits_{p\equiv1,7\operatorname*{mod}8}}
\left(  1-\frac1{p^{2}}\right)  ^{1/2}.
\end{align*}
As before,
\[
\Lambda_{-2}(s)=\left(  1-2^{-s}\right)  ^{-1}%
{\displaystyle\prod\limits_{\substack{p\equiv1,7 \\\operatorname*{mod}8 }}}
\left(  1-p^{-s}\right)  ^{-1}\cdot%
{\displaystyle\prod\limits_{\substack{p\equiv3,5 \\\operatorname*{mod}8 }}}
\left(  1-p^{-2s}\right)  ^{-1},
\]
\[
\Lambda_{-2}(s)^{2}=\zeta(s)L_{8}(s)\left(  1-2^{-s}\right)  ^{-1}%
{\displaystyle\prod\limits_{\substack{p\equiv3,5 \\\operatorname*{mod}8 }}}
\left(  1-p^{-2s}\right)  ^{-1}
\]
and, upon division, we obtain
\[
\underset{
\begin{array}
[c]{c}%
\downarrow\\
\sqrt{\pi}\kappa_{-2}%
\end{array}
}{\underbrace{\Lambda_{-2}(s)(s-1)^{1/2}}}\,=\,\underset{
\begin{array}
[c]{c}%
\downarrow\\
1
\end{array}
}{\underbrace{\zeta(s)(s-1)}}\,\underset{
\begin{array}
[c]{c}%
\downarrow\\
\frac{\ln\left(  1+\sqrt{2}\right)  }{\sqrt{2}}%
\end{array}
}{\underbrace{L_{8}(s)}}%
{\displaystyle\prod\limits_{\substack{p\equiv1,7 \\\operatorname*{mod}8 }}}
\left(  1-p^{-s}\right)  \cdot(s-1)^{-1/2}.
\]
Therefore
\[
\lim_{s\rightarrow1}%
{\displaystyle\prod\limits_{\substack{p\equiv1,3 \\\operatorname*{mod}8 }}}
\left(  1-\frac1{p^{s}}\right)  ^{-2}\cdot(s-1)=\frac\pi8\frac1{\kappa_{2}%
^{2}}=\frac{2^{1/2}}\pi%
{\displaystyle\prod\limits_{\substack{p\equiv1,3 \\\operatorname*{mod}8 }}}
\left(  1-\frac1{p^{2}}\right)  ^{-1},
\]
\[
\lim_{s\rightarrow1}%
{\displaystyle\prod\limits_{\substack{p\equiv1,7 \\\operatorname*{mod}8 }}}
\left(  1-\frac1{p^{s}}\right)  ^{-2}\cdot(s-1)=\frac{\ln\left(  1+\sqrt
{2}\right)  ^{2}}{2\pi}\frac1{\kappa_{-2}^{2}}=\frac{2^{3/2}\ln\left(
1+\sqrt{2}\right)  }{\pi^{2}}%
{\displaystyle\prod\limits_{\substack{p\equiv1,7 \\\operatorname*{mod}8 }}}
\left(  1-\frac1{p^{2}}\right)  ^{-1}
\]
but these alone do not go far enough. A slightly revised formula (\ref{E2}):
\[
\lim_{s\rightarrow1}%
{\displaystyle\prod\limits_{\substack{p\equiv1,5 \\\operatorname*{mod}8 }}}
\left(  1-\frac1{p^{s}}\right)  ^{-2}\cdot(s-1)=\frac\pi{16}\frac1{\kappa
_{1}^{2}}=\frac1\pi%
{\displaystyle\prod\limits_{\substack{p\equiv1,5 \\\operatorname*{mod}8 }}}
\left(  1-\frac1{p^{2}}\right)  ^{-1},
\]
when coupled with the preceding two limits, make possible the isolation of
$p\equiv1\operatorname*{mod}8$ as follows:
\begin{align*}
&  \ \frac4{\pi^{4}}\ln\left(  1+\sqrt{2}\right)  \overset{=\pi^{2}%
/8}{\overbrace{%
{\textstyle\prod\limits_{1,3,5,7}}
(1-p^{-2})^{-1}}}\cdot%
{\textstyle\prod\limits_{1}}
(1-p^{-2})^{-2}\\
\  &  =\frac{2^{1/2}}\pi%
{\textstyle\prod\limits_{1,3}}
(1-p^{-2})^{-1}\cdot\frac1\pi%
{\textstyle\prod\limits_{1,5}}
(1-p^{-2})^{-1}\cdot\frac{2^{3/2}\ln\left(  1+\sqrt{2}\right)  }{\pi^{2}}%
{\textstyle\prod\limits_{1,7}}
(1-p^{-2})^{-1}\\
\  &  =\lim_{s\rightarrow1}\left(
{\textstyle\prod\limits_{1,3}}
(1-p^{-s})^{-2}\cdot(s-1)\right)  \left(
{\textstyle\prod\limits_{1,5}}
(1-p^{-s})^{-2}\cdot(s-1)\right)  \left(
{\textstyle\prod\limits_{1,7}}
(1-p^{-s})^{-2}\cdot(s-1)\right) \\
\  &  =\,\underset{=1/4}{\underbrace{\lim_{s\rightarrow1}\left(
{\textstyle\prod\limits_{1,3,5,7}}
(1-p^{-s})^{-2}\cdot(s-1)^{2}\right)  }}\cdot\lim_{s\rightarrow1}\left(
{\textstyle\prod\limits_{1}}
(1-p^{-s})^{-4}\cdot(s-1)\right)  .
\end{align*}
Hence formula (\ref{E4}) is true.

\section{Unanswered Questions}

For the cases $\ell=5,6,7$, the function $\zeta(s)$ appears with exponents
$1,3,1$ respectively and \cite{SS}
\[
\kappa_{-3}=\tfrac12\sqrt{\tfrac{\sqrt{3}\ln\left(  2+\sqrt{3}\right)  }\pi}%
{\textstyle\prod\limits_{\substack{p\equiv5,7 \\\operatorname*{mod}12 }}}
\left(  1-\tfrac1{p^{2}}\right)  ^{-1/2},
\]
\[%
\begin{array}
[c]{lll}%
\kappa_{5}=\frac12\sqrt{\frac{\sqrt{5}}2}%
{\textstyle\prod\limits_{\substack{p\equiv11,13,17,19 \\\operatorname*{mod}20
}}}
\left(  1-\frac1{p^{2}}\right)  ^{-1/2}, &  & \kappa_{-5}=\frac12\sqrt
{\frac{\sqrt{5}\ln\left(  9+4\sqrt{5}\right)  }{3\pi}}%
{\textstyle\prod\limits_{\substack{p\equiv2,3 \\\operatorname*{mod}5 }}}
\left(  1-\frac1{p^{2}}\right)  ^{-1/2},\\
\kappa_{6}=\frac12\sqrt{\frac{\sqrt{6}}2}%
{\textstyle\prod\limits_{\substack{p\equiv13,17,19,23 \\\operatorname*{mod}24
}}}
\left(  1-\frac1{p^{2}}\right)  ^{-1/2}, &  & \kappa_{-6}=\frac12\sqrt
{\frac{\sqrt{6}\ln\left(  5+2\sqrt{6}\right)  }{2\pi}}%
{\textstyle\prod\limits_{\substack{p\equiv7,11,13,17 \\\operatorname*{mod}24
}}}
\left(  1-\frac1{p^{2}}\right)  ^{-1/2},\\
\kappa_{7}=\frac34\sqrt{\frac{\sqrt{7}}6}%
{\textstyle\prod\limits_{\substack{p\equiv3,5,6 \\\operatorname*{mod}7 }}}
\left(  1-\frac1{p^{2}}\right)  ^{-1/2}, &  & \kappa_{-7}=\frac12\sqrt
{\frac{\sqrt{7}\ln\left(  8+3\sqrt{7}\right)  }{3\pi}}%
{\textstyle\prod\limits_{\substack{p\equiv5,11,13, \\15,17,23
\\\operatorname*{mod}28 }}}
\left(  1-\frac1{p^{2}}\right)  ^{-1/2}.
\end{array}
\]
It does not seem to be possible to isolate $p\equiv1\operatorname*{mod}\ell$
beyond the following partial results:
\[
\lim_{s\rightarrow1}%
{\displaystyle\prod\limits_{\substack{p\equiv1,4 \\\operatorname*{mod}5 }}}
\left(  1-\frac1{p^{s}}\right)  ^{-2}\cdot(s-1)=\frac{\ln\left(  9+4\sqrt
{5}\right)  ^{2}}{45\pi}\frac1{\kappa_{-5}^{2}}=\frac{\sqrt{5}\ln\left(
9+4\sqrt{5}\right)  }{3\pi^{2}}%
{\displaystyle\prod\limits_{\substack{p\equiv1,4 \\\operatorname*{mod}5 }}}
\left(  1-\frac1{p^{2}}\right)  ^{-1},
\]
\[
\lim_{s\rightarrow1}%
{\displaystyle\prod\limits_{\substack{p\equiv1,2,4 \\\operatorname*{mod}7 }}}
\left(  1-\frac1{p^{s}}\right)  ^{-2}\cdot(s-1)=\frac{9\pi}{112\pi}%
\frac1{\kappa_{7}^{2}}=\frac{3\sqrt{7}}{4\pi}%
{\displaystyle\prod\limits_{\substack{p\equiv1,2,4 \\\operatorname*{mod}7 }}}
\left(  1-\frac1{p^{2}}\right)  ^{-1}
\]
for $\ell=5,7$, which are deduced from
\[
\Lambda_{-5}(s)=\left(  1-5^{-s}\right)  ^{-1}%
{\displaystyle\prod\limits_{\substack{p\equiv1,4 \\\operatorname*{mod}5 }}}
\left(  1-p^{-s}\right)  ^{-1}\cdot%
{\displaystyle\prod\limits_{\substack{p\equiv2,3 \\\operatorname*{mod}5 }}}
\left(  1-p^{-2s}\right)  ^{-1},
\]
\[
\Lambda_{-5}(s)^{2}=\zeta(s)L_{20}(s)\left(  1-5^{-s}\right)  ^{-1}%
{\displaystyle\prod\limits_{\substack{p\equiv2,3 \\\operatorname*{mod}5 }}}
\left(  1-p^{-2s}\right)  ^{-1}\cdot\left(  1+2^{-s}\right)  ^{-1}
\]
and
\[
\Lambda_{7}(s)=\left(  1-7^{-s}\right)  ^{-1}%
{\displaystyle\prod\limits_{\substack{2<p\equiv1,2,4 \\\operatorname*{mod}7
}}}
\left(  1-p^{-s}\right)  ^{-1}\cdot%
{\displaystyle\prod\limits_{\substack{p\equiv3,5,6 \\\operatorname*{mod}7 }}}
\left(  1-p^{-2s}\right)  ^{-1}\cdot\left[  \left(  1-2^{-s}\right)
^{-1}-2^{-s}\right]  ,
\]
\[
\Lambda_{7}(s)^{2}=\zeta(s)L_{-28}(s)\left(  1-7^{-s}\right)  ^{-1}\cdot%
{\displaystyle\prod\limits_{\substack{p\equiv3,5,6 \\\operatorname*{mod}7 }}}
\left(  1-p^{-2s}\right)  ^{-1}\cdot\left[  \left(  1-2^{-s}\right)
^{-1}-2^{-s}\right]  ^{2}\left(  1-2^{-s}\right)
\]
as before. The L-series $\Lambda_{-3}(s),\Lambda_{5}(s),\Lambda_{6}%
(s),\Lambda_{-6}(s),\Lambda_{-7}(s)$ might be difficult to study, due to a
failure of $f(n)$-multiplicativity. When $m=6$, for example, $f(10)=1$ since
$10=2^{2}+6\cdot1^{2}$, but $f(2)=f(5)=0$. When $m=-6$, as another example,
$f(10)=1$ since $10=4^{2}-6\cdot1^{2}$, but again $f(2)=f(5)=0$. Although the
asymptotics of $%
{\textstyle\sum\nolimits_{n\leq N}}
b(n)$ for $\ell=5,6,7$ are understood, expressions for leading coefficients
(analogous to those for $\ell=3,4,8$) remain open.

\section{Acknowledgement}

I\ thank Pascal Sebah, my coauthor in \cite{FS}, for his skillful numerical
computations over many years! There is a way to avoid the \textquotedblleft
roundabout\textquotedblright\ calculation of residues (involving \cite{SS}, as
presented here) and we revisit this topic in \cite{FS2}. \ Our work is
extended in \cite{FMS} and we gratefully acknowledge Greg Martin for his
mastery of the subject.

\end{document}